\documentclass{article}
\usepackage{graphicx}
\usepackage{amsmath}
\usepackage{amsfonts}
\usepackage{amssymb}
\newtheorem{theorem}{Theorem}

\newtheorem{corollary}[theorem]{Corollary}

\newtheorem{proposition}[theorem]{Proposition}

\newenvironment{proof}[1][Proof]{\textbf{#1.} }{\ \rule{0.5em}{0.5em}}

\begin{document}

\title{Covariant Poisson Structures on Complex Projective Spaces}
\author{Albert Jeu-Liang Sheu\thanks{Partially supported by NSF Grant DMS-9623008.}\\Department of Mathematics\\University of Kansas\\Lawrence, KS 66045\\U. S. A.}
\date{}
\maketitle

\section{Introduction}

One of the most intriguing aspects of the theory of quantum groups and quantum
spaces \cite{D,FRT,So,Wo:ts,Wo:cm,Po,Ri:ncq,VaSo:af} is the close interplay
between the geometric structure of the underlying Poisson Lie group (or
Poisson space) \cite{We,LuWe:plg} and the algebraic structure on the
corresponding quantum group (or quantum space). For example, Soibelman's
classification \cite{So} of all irreducible *-representations of the quantum
algebra $C\left(  G_{q}\right)  $ for complex Poisson simple Lie groups $G$
gives a one-to-one correspondence between irreducible *-representations of
$C\left(  G_{q}\right)  $ and the symplectic leaves \cite{We} on $G$. This
leads to a groupoid C*-algebraic \cite{Re} approach to study the structure of
the algebra $C\left(  G_{q}\right)  $ \cite{Sh:cqg} which shows that the
decomposition of $SU\left(  n\right)  $ (or $\mathbb{S}^{2n+1}$) by symplectic
leaves of various dimensions corresponds to a compatible decomposition of
$C\left(  SU\left(  n\right)  _{q}\right)  $ (or $C\left(  \mathbb{S}%
_{q}^{2n+1}\right)  $) by its (closed) ideals in the spirit of noncommutative geometry.

It is well known that the standard multiplicative Poisson structure on
$SU\left(  n+1\right)  $ induces a standard covariant Poisson structure on the
homogeneous spaces $\mathbb{S}^{2n+1}=SU\left(  n+1\right)  /SU\left(
n\right)  $ and $\mathbb{C}P^{n}=SU\left(  n+1\right)  /U\left(  n\right)  $
determined by the Poisson Lie subgroups $SU\left(  n\right)  $ and $U\left(
n\right)  $, respectively. On the other hand, Lu and Weinstein \cite{LuWe:ps}
described explicitly all $SU\left(  2\right)  $-covariant Poisson structures
on $\mathbb{S}^{2}=\mathbb{C}P^{1}$ including a one-parameter family of
nonstandard $SU\left(  2\right)  $-covariant Poisson structures on
$\mathbb{S}^{2}$, and showed that each nonstandard covariant Poisson sphere
$\mathbb{S}_{c}^{2}$ contains a copy of the trivial Poisson $1$-sphere
$\mathbb{S}^{1}$ (consisting of a circle family of $0$-dimensional symplectic
leaves) and exactly two $2$-dimensional symplectic leaves. This geometric
structure is reflected faithfully in the algebraic structure of the algebra
$C\left(  \mathbb{S}_{qc}^{2}\right)  $ of the nonstandard quantum spheres
$\mathbb{S}_{qc}^{2}$ \cite{Sh:qp}.

Dijkhuizen and Noumi studied in great detail \cite{DiNo} a one-parameter
family of nonstandard quantum projective spaces $\mathbb{C}P_{q,c}^{n}$ with
quantum algebras $C\left(  \mathbb{C}P_{q,c}^{n}\right)  $. In \cite{Sh:qcp},
the structure of $C\left(  \mathbb{C}P_{q,c}^{n}\right)  $ is studied and
analyzed as a groupoid C*-algebra, and an algebraic decomposition of $C\left(
\mathbb{C}P_{q,c}^{n}\right)  $ by a closed ideal indicates that the
underlying nonstandard Poisson $\mathbb{C}P^{n}$ should contain an embeded
copy of the standard Poisson $\mathbb{S}^{2n-1}$. In this paper, we describe a
one-parameter family of nonstandard $SU\left(  n+1\right)  $-covariant Poisson
structures $\tau_{c}$ on the projective space $\mathbb{C}P^{n}$ that
represents the classical counterpart of the quantum family $C\left(
\mathbb{C}P_{q,c}^{n}\right)  $, and show that the standard Poisson
$\mathbb{S}^{2n-1}$ is indeed embedded in each of these nonstandard Poisson
$\mathbb{C}P^{n}$. We also show that (non-zero) $SU\left(  n\right)
$-invaraint contravariant 2-tensors on $\mathbb{S}^{2n-1}$ with $n\neq3$ (or
$\mathbb{C}P^{n}$) are unique, up to a constant factor. We remark that in
\cite{KRR} Khoroshkin, Radul, and Rubtsov obtained interesting results about
covariant Poisson structures on coadjoint orbits, including $\mathbb{C}P^{n}$.
Our approach, motivated by Dijkhuizen and Noumi's work \cite{DiNo}, is
different from theirs and the embedding of the standard Poisson $\mathbb{S}%
^{2n-1}$ in the nonstandard Poisson $\mathbb{C}P^{n}$ is new.

\section{Poisson structure on Lie groups}

In this section, we discuss some basic properties of affine Poisson structures
in the form needed later. We recall that an affine Poisson structure on a Lie
group $G$ is given by a Poisson 2-tensor $\pi\in\Gamma\left(  \wedge
^{2}TG\right)  $, such that
\[
\pi\left(  gh\right)  =L_{g}\left(  \pi\left(  h\right)  \right)
+R_{h}\left(  \pi\left(  g\right)  \right)  -L_{g}R_{h}\left(  \pi\left(
e\right)  \right)
\]
for any $g,h\in G$ \cite{We:aps}, or equivalently,
\[
\pi_{l}\left(  g\right)  :=\pi\left(  g\right)  -L_{g}\left(  \pi\left(
e\right)  \right)
\]
for $g\in G$ defines a multiplicative Poisson 2-tensor on $G$
\cite{Lu:maps,DaSo}, where $L_{g}$ and $R_{g}$ are the left and the right
actions by $g\in G$, respectively, and $e$ is the identity element of $G$. For
an affine Poisson 2-tensor $\pi$ on a Lie group $G$, the left action of the
Poisson-Lie group $\left(  G,\pi_{l}\right)  $ on the Poisson manifold
$\left(  G,\pi\right)  $ by left translation is a Poisson action, i.e. the
multiplication map $G\times G\rightarrow G$ is a Poisson map, where $G\times
G$ and $G$ are endowed with the product Poisson structures $\pi_{l}\times\pi$
and $\pi$, respectively. In another word, $\pi$ on $G$ (as a homogeneous space
of $G$) is a (left) $\left(  G,\pi_{l}\right)  $-covariant Poisson structure.

A typical example of an affine Poisson structure on a Poisson-Lie group $G$
with multiplicative Poisson 2-tensor $\pi$ is provided by a right translation
$\pi_{\sigma}$ of $\pi$ by an elemnet $\sigma\in G$, i.e.
\[
\pi_{\sigma}\left(  g\right)  :=R_{\sigma}\left(  \pi\left(  g\sigma
^{-1}\right)  \right)
\]
for $g\in G$. Since the right translation by $\sigma$ on $G$ is a
diffeomorphism on $G$, the `push-forward' $\pi_{\sigma}$ of $\pi$ by
$R_{\sigma}$ is clearly a Poisson 2-tensor on $G$. Furthermore,
\[
\left(  \pi_{\sigma}\right)  _{l}\left(  g\right)  =R_{\sigma}\left(
\pi\left(  g\sigma^{-1}\right)  \right)  -L_{g}\left(  R_{\sigma}\left(
\pi\left(  \sigma^{-1}\right)  \right)  \right)
\]%
\[
=R_{\sigma}\left(  L_{g}\left(  \pi\left(  \sigma^{-1}\right)  \right)
+R_{\sigma^{-1}}\left(  \pi\left(  g\right)  \right)  \right)  -L_{g}\left(
R_{\sigma}\left(  \pi\left(  \sigma^{-1}\right)  \right)  \right)
\]%
\[
=R_{\sigma}L_{g}\left(  \pi\left(  \sigma^{-1}\right)  \right)  +\pi\left(
g\right)  -R_{\sigma}L_{g}\left(  \pi\left(  \sigma^{-1}\right)  \right)
=\pi\left(  g\right)
\]
which is a multiplicative Poisson 2-tensor on $G$. So $\pi_{\sigma}$ on $G$
(as a homogeneous space of $G$) is a (left) $\left(  G,\pi\right)  $-covariant
Poisson structure, for any $\sigma\in G$. Note that
\[
\pi_{\sigma}=\pi+\left(  X_{\sigma}\right)  ^{l}%
\]
where $X_{\sigma}:=\pi_{\sigma}\left(  e\right)  =R_{\sigma}\left(  \pi\left(
\sigma^{-1}\right)  \right)  \in\frak{g}\wedge\frak{g}$ and $X^{l}\left(
g\right)  :=L_{g}\left(  X\right)  $ is the left-invariant 2-tensor generated
by $X\in\frak{g}\wedge\frak{g}$, since $\left(  \pi_{\sigma}\right)  _{l}=\pi$.

When a closed subgroup $H$ of a Lie group $G\ $ is coisotropic \cite{We:sm}
with respect to a Poisson structure $\rho$ on $G$, i.e.
\[
\rho\left(  gh\right)  -R_{h}\left(  \rho\left(  g\right)  \right)  \in
L_{gh}\left(  \frak{h}\wedge\frak{g}\right)
\]
for all $g\in G$ and $h\in H$, it is easy to see that the Poisson bracket
$\left\{  f_{1},f_{2}\right\}  :=\left(  df_{1}\wedge df_{2}\right)  \left(
\rho\right)  $ of $f_{1},f_{2}\in C^{\infty}\left(  G/H\right)  \subset
C^{\infty}\left(  G\right)  $ is still in $C^{\infty}\left(  G/H\right)  $ and
hence induces a Poisson structure on $G/H$, or equivalently, a Poisson
2-tensor $\tilde{\rho}$ on the homogeneous space $G/H$ is well defined by
\[
\tilde{\rho}\left(  \left[  gH\right]  \right)  :=\left[  \rho\left(
g\right)  \right]  \in L_{g}\left(  \wedge^{2}\left(  \frak{g}/\frak{h}%
\right)  \right)  =\wedge^{2}T_{\left[  gH\right]  }\left(  G/H\right)  .
\]

Given a Poisson-Lie group $\left(  G,\pi\right)  $ and $\sigma\in G$, if a
closed subgroup $H$ of $G$ is coisotropic with respect to $\pi$, then $H$ is
coisotropic with respect to the affine Poisson structure $\pi_{\sigma}$ on $G$
if and only if
\[
\left(  X_{\sigma}\right)  ^{l}\left(  gh\right)  -R_{h}\left(  \left(
X_{\sigma}\right)  ^{l}\left(  g\right)  \right)  \in L_{gh}\left(
\frak{h}\wedge\frak{g}\right)  ,
\]
since $\pi_{\sigma}=\pi+\left(  X_{\sigma}\right)  ^{l}$ and $\pi\left(
gh\right)  -R_{h}\left(  \pi\left(  g\right)  \right)  \in L_{gh}\left(
\frak{h}\wedge\frak{g}\right)  $. Now
\[
\left(  X_{\sigma}\right)  ^{l}\left(  gh\right)  -R_{h}\left(  \left(
X_{\sigma}\right)  ^{l}\left(  g\right)  \right)  =L_{gh}\left(  X_{\sigma
}\right)  -R_{h}\left(  L_{g}\left(  X_{\sigma}\right)  \right)
\]%
\[
=L_{g}\left(  L_{h}X_{\sigma}-R_{h}X_{\sigma}\right)  =L_{g}L_{h}\left(
X_{\sigma}-L_{h^{-1}}R_{h}X_{\sigma}\right)
\]%
\[
=L_{gh}\left(  \operatorname{id}-\operatorname{Ad}_{h^{-1}}\right)  \left(
X_{\sigma}\right)  .
\]
So $\left(  X_{\sigma}\right)  ^{l}\left(  gh\right)  -R_{h}\left(  \left(
X_{\sigma}\right)  ^{l}\left(  g\right)  \right)  \in L_{gh}\left(
\frak{h}\wedge\frak{g}\right)  $ for all $\left(  g,h\right)  \in G\times H$
if and only if
\[
\left(  \operatorname{id}-\operatorname{Ad}_{h^{-1}}\right)  \left(
X_{\sigma}\right)  \in\frak{h}\wedge\frak{g}%
\]
for all $\left(  g,h\right)  \in G\times H$, or equivalently,
\[
\operatorname{ad}_{\frak{h}}\left(  X_{\sigma}\right)  \subset\frak{h}%
\wedge\frak{g},
\]
since $\operatorname{ad}_{\frak{h}}\left(  \frak{h}\wedge\frak{g}\right)
\subset\frak{h}\wedge\frak{g}$ and $\operatorname{Ad}_{H}\left(
\frak{h}\wedge\frak{g}\right)  \subset\frak{h}\wedge\frak{g}$. Thus we get the
following result.

\begin{proposition}
Given a Poisson-Lie group $\left(  G,\pi\right)  $ and a closed subgroup $H$
of $G$ that is coisotropic with respect to $\pi$, the subgroup $H$ is
coisotropic with respect to $\pi_{\sigma}$ for $\sigma\in G$, if and only if
$\operatorname{ad}_{\frak{h}}\left(  X_{\sigma}\right)  \subset\frak{h}%
\wedge\frak{g}$, where $X_{\sigma}:=R_{\sigma}\left(  \pi\left(  \sigma
^{-1}\right)  \right)  \in\frak{g}\wedge\frak{g}$.
\end{proposition}

In case the multiplicative Poisson structure $\pi$ on $G$ is given by an
$r$-matrix $r\in\frak{g}\wedge\frak{g}$ (satisfying the modified Yang-Baxter
equation), i.e.
\[
\pi\left(  g\right)  =L_{g}r-R_{g}r,
\]
we have
\[
X_{\sigma}=R_{\sigma}\left(  \pi\left(  \sigma^{-1}\right)  \right)
=R_{\sigma}\left(  L_{\sigma^{-1}}r-R_{\sigma^{-1}}r\right)
=\operatorname{Ad}_{\sigma^{-1}}\left(  r\right)  -r.
\]
A closed subgroup $H$ being coisotropic with respect to $\pi$ is equivalent
to
\[
\operatorname{ad}_{\frak{h}}\left(  r\right)  \subset\frak{h}\wedge\frak{g},
\]
since
\[
\pi\left(  gh\right)  -R_{h}\left(  \pi\left(  g\right)  \right)
=L_{gh}r-R_{gh}r-R_{h}\left(  L_{g}r-R_{g}r\right)
\]%
\[
=L_{gh}r-R_{gh}r-L_{g}R_{h}r+R_{h}R_{g}r=L_{gh}r-L_{g}R_{h}r
\]%
\[
=L_{gh}\left(  r-L_{h^{-1}}R_{h}r\right)  =L_{gh}\left(  \operatorname{id}%
-\operatorname{Ad}_{h^{-1}}\right)  \left(  r\right)
\]
and $L_{gh}\left(  \operatorname{id}-\operatorname{Ad}_{h^{-1}}\right)
\left(  r\right)  \in L_{gh}\left(  \frak{h}\wedge\frak{g}\right)  $ for all
$\left(  g,h\right)  \in G\times H$ if and only if $\left(  \operatorname{id}%
-\operatorname{Ad}_{h^{-1}}\right)  \left(  r\right)  \in\frak{h}%
\wedge\frak{g}$ for all $h\in H$.

\begin{corollary}
Given a Poisson-Lie group $\left(  G,\pi\right)  $ with $\pi$ defined by an
$r$-matrix $r\in\frak{g}\wedge\frak{g}$, a $\pi$-coisotropic closed subgroup
$H$ of $G$ is coisotropic with respect to $\pi_{\sigma}$ for $\sigma\in G$, if
and only if $\operatorname{ad}_{\frak{h}}\left(  \operatorname{Ad}%
_{\sigma^{-1}}\left(  r\right)  \right)  \subset\frak{h}\wedge\frak{g}$.
\end{corollary}

\section{Non-standard Poisson $\mathbb{C}P^{n}$}

Recall that the satandard Poisson $SU\left(  n\right)  $ is defined (up to a
constant multiple) by the Poisson 2-tensor $\pi\left(  u\right)  =\pi^{\left(
n\right)  }\left(  u\right)  :=L_{u}r-R_{u}r$ determined by the $r$-matrix
\[
r:=\sum_{1\leq i<j\leq n}X_{ij}^{+}\wedge X_{ij}^{-}%
\]
where $X_{ij}^{+}=e_{ij}-e_{ji}$, $X_{ij}^{-}=i\left(  e_{ij}+e_{ji}\right)
$, and $e_{ij}$ are the matrix units.

It is well known that $SU\left(  n-1\right)  =\left\{  1\right\}  \oplus
SU\left(  n-1\right)  $ (or $SU\left(  n-1\right)  \oplus\left\{  1\right\}
$) and
\[
U\left(  n-1\right)  :=\left\{  \det\left(  u\right)  ^{-1}\oplus u:u\in
SU\left(  n-1\right)  \right\}
\]
are Poisson-Lie subgroups of $SU\left(  n\right)  $ and hence induce the
`standard' $SU\left(  n\right)  $-covariant Poisson structures $\rho
=\rho^{\left(  n\right)  }$ and $\tau=\tau^{\left(  n-1\right)  }$ on the
sphere
\[
\mathbb{S}^{2n-1}\cong SU\left(  n\right)  /\left[  \left\{  1\right\}  \oplus
SU\left(  n-1\right)  \right]
\]
and the complex projective space
\[
\mathbb{C}P^{n-1}\cong SU\left(  n\right)  /U\left(  n-1\right)  ,
\]
respectively.

\begin{theorem}
The closed subgroup $U\left(  n-1\right)  $ of $SU\left(  n\right)  $ is
coisotropic with respect to the (left) $SU\left(  n\right)  $-covariant affine
Poisson structure $\pi_{\sigma_{c}}$ on $SU\left(  n\right)  $ defined by
\[
\sigma_{c}:=\left(  \sqrt{c}e_{11}+\sqrt{1-c}e_{n1}-\sqrt{1-c}e_{1n}+\sqrt
{c}e_{nn}\right)  +\sum_{i=2}^{n-1}e_{ii}\in SU\left(  n\right)
\]
with $c\in\left[  0,1\right]  $. Hence $\pi_{\sigma_{c}}$ induces a (left)
$SU\left(  n\right)  $-covariant Poisson structure $\tau_{c}$ on
$\mathbb{C}P^{n-1}\cong SU\left(  n\right)  /U\left(  n-1\right)  $.
\end{theorem}

\begin{proof}
We set $\sigma=\sigma_{c}$ for simplicity. It is easy to see that if the
Poisson structure $\tau_{c}$ induced by $\pi_{\sigma_{c}}$ on $\mathbb{C}%
P^{n-1}\cong SU\left(  n\right)  /U\left(  n-1\right)  $ is well defined, then
it is automatically (left) $SU\left(  n\right)  $-covariant since $\pi
_{\sigma_{c}}$ is.

Now since $U\left(  n-1\right)  $ is a Poisson-Lie subgroup and hence
coisotropic with respect to $\pi$, we have
\[
\operatorname{ad}_{\frak{u}\left(  n-1\right)  }\left(  r\right)
\subset\frak{u}\left(  n-1\right)  \wedge\frak{su}\left(  n\right)  ,
\]
and hence only need to show that
\[
\operatorname{ad}_{\frak{u}\left(  n-1\right)  }\left(  \operatorname{Ad}%
_{\sigma^{-1}}\left(  r\right)  \right)  \subset\frak{u}\left(  n-1\right)
\wedge\frak{su}\left(  n\right)  .
\]

From
\[
\left\{
\begin{array}
[c]{ll}%
\operatorname{Ad}_{\sigma^{-1}}\left(  X_{ij}^{+}\right)  =X_{ij}^{+}, &
\text{if }1<i<j<n\\
\operatorname{Ad}_{\sigma^{-1}}\left(  X_{1j}^{+}\right)  =\sqrt{c}X_{1j}%
^{+}+\sqrt{1-c}X_{jn}^{+}, & \text{if }1<j<n\\
\operatorname{Ad}_{\sigma^{-1}}\left(  X_{in}^{+}\right)  =-\sqrt{1-c}%
X_{1i}^{+}+\sqrt{c}X_{in}^{+}, & \text{if }1<i<n\\
\operatorname{Ad}_{\sigma^{-1}}\left(  X_{1n}^{+}\right)  =X_{1n}^{+}, &
\end{array}
\right.
\]
and
\[
\left\{
\begin{array}
[c]{ll}%
\operatorname{Ad}_{\sigma^{-1}}\left(  X_{ij}^{-}\right)  =X_{ij}^{-}, &
\text{if }1<i<j<n\\
\operatorname{Ad}_{\sigma^{-1}}\left(  X_{1j}^{-}\right)  =\sqrt{c}X_{1j}%
^{-}-\sqrt{1-c}X_{jn}^{-}, & \text{if }1<j<n\\
\operatorname{Ad}_{\sigma^{-1}}\left(  X_{in}^{-}\right)  =\sqrt{1-c}%
X_{1i}^{-}+\sqrt{c}X_{in}^{-}, & \text{if }1<i<n\\
\operatorname{Ad}_{\sigma^{-1}}\left(  X_{1n}^{-}\right)  =\left(
2c-1\right)  X_{1n}^{-}+2i\sqrt{c\left(  1-c\right)  }\left(  e_{11}%
-e_{nn}\right)  . &
\end{array}
\right.
\]
we get
\[
\operatorname{Ad}_{\sigma^{-1}}\left(  r\right)  =\sum_{1\leq i<j\leq
n}\operatorname{Ad}_{\sigma^{-1}}\left(  X_{ij}^{+}\right)  \wedge
\operatorname{Ad}_{\sigma^{-1}}\left(  X_{ij}^{-}\right)
\]%
\[
=\sum_{1<i<j<n}X_{ij}^{+}\wedge X_{ij}^{-}+X_{1n}^{+}\wedge\left[  \left(
2c-1\right)  X_{1n}^{-}+2\sqrt{c\left(  1-c\right)  }i\left(  e_{11}%
-e_{nn}\right)  \right]
\]%
\[
+\sum_{1<i<n}\left(  \sqrt{c}X_{1i}^{+}+\sqrt{1-c}X_{in}^{+}\right)
\wedge\left(  \sqrt{c}X_{1i}^{-}-\sqrt{1-c}X_{in}^{-}\right)
\]%
\[
+\sum_{1<i<n}\left(  -\sqrt{1-c}X_{1i}^{+}+\sqrt{c}X_{in}^{+}\right)
\wedge\left(  \sqrt{1-c}X_{1i}^{-}+\sqrt{c}X_{in}^{-}\right)
\]%
\[
=\sum_{1<i<j<n}X_{ij}^{+}\wedge X_{ij}^{-}+2\sqrt{c\left(  1-c\right)  }%
X_{1n}^{+}\wedge i\left(  e_{11}-e_{nn}\right)
\]%
\[
+\left(  2c-1\right)  \left[  \left(  X_{1n}^{+}\wedge X_{1n}^{-}\right)
+\sum_{1<i<n}\left(  X_{1i}^{+}\wedge X_{1i}^{-}+X_{in}^{+}\wedge X_{in}%
^{-}\right)  \right]
\]%
\[
+2\sqrt{c\left(  1-c\right)  }\sum_{1<i<n}\left(  X_{in}^{+}\wedge X_{1i}%
^{-}-X_{1i}^{+}\wedge X_{in}^{-}\right)
\]%
\[
=\left(  1-\left(  2c-1\right)  \right)  \sum_{1<i<j<n}X_{ij}^{+}\wedge
X_{ij}^{-}+2\sqrt{c\left(  1-c\right)  }X_{1n}^{+}\wedge i\left(
e_{11}-e_{nn}\right)
\]%
\[
+\left(  2c-1\right)  \sum_{1\leq i<j\leq n}X_{ij}^{+}\wedge X_{ij}^{-}%
+2\sqrt{c\left(  1-c\right)  }\sum_{1<i<n}\left(  X_{in}^{+}\wedge X_{1i}%
^{-}-X_{1i}^{+}\wedge X_{in}^{-}\right)
\]%
\[
=2\left(  1-c\right)  \sum_{1<i<j<n}X_{ij}^{+}\wedge X_{ij}^{-}+2\sqrt
{c\left(  1-c\right)  }X_{1n}^{+}\wedge i\left(  e_{11}-e_{nn}\right)
\]%
\[
+\left(  2c-1\right)  r+2\sqrt{c\left(  1-c\right)  }\sum_{1<i<n}\left(
X_{in}^{+}\wedge X_{1i}^{-}-X_{1i}^{+}\wedge X_{in}^{-}\right)
\]
which is in $\left(  2c-1\right)  r+\left(  \frak{u}\left(  n-1\right)
\wedge\frak{su}\left(  n\right)  \right)  $, since $X_{ij}^{+},X_{ij}%
^{-},X_{in}^{+},X_{in}^{-},i\left(  e_{11}-e_{nn}\right)  \in\frak{u}\left(
n-1\right)  $ for all $1<i<j<n$. So we get
\[
\operatorname{ad}_{\frak{u}\left(  n-1\right)  }\left(  \operatorname{Ad}%
_{\sigma^{-1}}\left(  r\right)  \right)  \subset\left(  2c-1\right)
\operatorname{ad}_{\frak{u}\left(  n-1\right)  }\left(  r\right)
+\operatorname{ad}_{\frak{u}\left(  n-1\right)  }\left(  \frak{u}\left(
n-1\right)  \wedge\frak{su}\left(  n\right)  \right)
\]%
\[
\subset\frak{u}\left(  n-1\right)  \wedge\frak{su}\left(  n\right)  ,
\]
because $\operatorname{ad}_{\frak{u}\left(  n-1\right)  }\left(  r\right)
\subset\frak{u}\left(  n-1\right)  \wedge\frak{su}\left(  n\right)  $ and
\[
\operatorname{ad}_{\frak{u}\left(  n-1\right)  }\left(  \frak{u}\left(
n-1\right)  \wedge\frak{su}\left(  n\right)  \right)  \subset\frak{u}\left(
n-1\right)  \wedge\frak{su}\left(  n\right)  .
\]
\end{proof}

The Poisson manifold $\left(  \mathbb{C}P^{n-1},\tau_{c}\right)  $ with
$c\in\left(  0,1\right)  $ is referred to as a nonstandard Poisson
$\mathbb{C}P^{n-1}$. Note that $\tau_{1}=\tau^{\left(  n-1\right)  }$ the
standard Poisson structure on $\mathbb{C}P^{n-1}$ since $\sigma_{1}=1\in
SU\left(  n\right)  $ and hence $\pi_{\sigma_{1}}=\pi$. On the other hand,
$R_{\sigma_{0}}$ simply swaps the first column with the $n$-th column and
hence $\tau_{0}$ is the standard Poisson structure on $\mathbb{C}P^{n-1}\cong
SU\left(  n\right)  /\left(  SU\left(  n-1\right)  \oplus\left\{  1\right\}
\right)  $ induced by $\pi$.

We remark that $X_{ij}^{+},X_{ij}^{-},X_{in}^{+},X_{in}^{-}\in\frak{su}\left(
n-1\right)  $ but $i\left(  e_{11}-e_{nn}\right)  \notin\frak{su}\left(
n-1\right)  $, so $\operatorname{ad}_{\frak{su}\left(  n-1\right)  }\left(
\operatorname{Ad}_{\sigma^{-1}}\left(  r\right)  \right)  \nsubseteqq
\frak{su}\left(  n-1\right)  \wedge\frak{su}\left(  n\right)  $ and hence
$\pi_{\sigma_{c}}$ does not induce a Poisson structure on $\mathbb{S}%
^{2n-1}\cong SU\left(  n\right)  /SU\left(  n-1\right)  $. On the other hand,
as a generalization of Lu and Weinstein's result on covariant Poisson spheres
$\mathbb{S}^{2}=\mathbb{C}P^{1}$ \cite{LuWe:ps}, we can show that $\left(
\mathbb{C}P^{n-1},\tau_{c}\right)  $ contains a copy of the standard Poisson
sphere $\left(  \mathbb{S}^{2n-3},\rho^{\left(  n-1\right)  }\right)  $. Here
it is understood that $\rho^{\left(  1\right)  }=0$ on $\mathbb{S}^{1}$ by definition.

\begin{theorem}
The standard Poisson sphere $\left(  \mathbb{S}^{2n-3},\rho^{\left(
n-1\right)  }\right)  $ is embedded in $\left(  \mathbb{C}P^{n-1},\tau
_{c}\right)  $ for $c\in\left(  0,1\right)  $ and $n\geq2$.
\end{theorem}

\begin{proof}
Note that the quotient map $\phi:SU\left(  n\right)  \rightarrow
\mathbb{C}P^{n-1}$ can be viewed as the composition of the quotient map
\[
\phi_{1}:u\in SU\left(  n\right)  \mapsto u_{1}\in\mathbb{S}^{2n-1}\cong
SU\left(  n\right)  /SU\left(  n-1\right)
\]
and the quotient map
\[
\phi_{2}:v\in\mathbb{S}^{2n-1}\mapsto\left[  v\right]  \in\mathbb{C}%
P^{n-1}\cong\mathbb{S}^{2n-1}/\mathbb{T},
\]
where the circle group $\mathbb{T}$ acts diagonally on $\mathbb{S}%
^{2n-1}\subset\mathbb{C}^{n}$ and
\[
u_{1}:=\left(  u_{11},u_{21},...,u_{n1}\right)  \in\mathbb{S}^{2n-1}%
\subset\mathbb{C}^{n}%
\]
is the first column of $u\in SU\left(  n\right)  $. It is well known that
$\phi_{2}$ is a diffeomorphism from the submanifold
\[
S_{+}:=\left\{  v\in\mathbb{S}^{2n-1}:v_{1}>0\right\}  \subset\mathbb{S}%
^{2n-1}%
\]
onto its image $\phi_{2}\left(  S_{+}\right)  $ $\subset\mathbb{C}P^{n-1}$,
and
\[
\phi_{3}:v\in S_{c}\mapsto\phi_{3}\left(  v\right)  :=\frac{1}{\sqrt{1-c}%
}\left(  v_{2},...,v_{n}\right)  \in\mathbb{S}^{2n-3}%
\]
is a diffeomorphism identifying
\[
S_{c}:=\left\{  v\in\mathbb{S}^{2n-1}:v_{1}=\sqrt{c}\right\}  \subset S_{+}%
\]
with $\mathbb{S}^{2n-3}$. We denote by $\psi:u\in SU\left(  n\right)  \mapsto
u_{n}\in\mathbb{S}^{2n-1}$ the projection to the last column. Functions
similar to $\phi_{1}$, $\phi_{2}$, and $\psi$, for other dimensions than $n$,
will be denoted by the same symbols for the simplicity of notation.

First we assume that $n>2$. For each $v\in S_{c}$, we can find some
$u^{\prime}\in SU\left(  n\right)  $ with the first column $u_{1}^{\prime
}=\left(  1,0,...,0\right)  $ and the last column $u_{n}^{\prime}=\sqrt
{1-c}^{-1}\left(  0,v_{2},...,v_{n}\right)  $. Note that the first row of
$u^{\prime}$ has to be $\left(  1,0,...,0\right)  $, and hence
\[
u^{\prime}=1\oplus u^{\prime\prime}\in\left\{  1\right\}  \oplus SU\left(
n-1\right)
\]
for some $u^{\prime\prime}\in SU\left(  n-1\right)  $ with
\[
\left(  u^{\prime\prime}\right)  _{n-1}=\sqrt{1-c}^{-1}\left(  v_{2}%
,...,v_{n}\right)  =\phi_{3}\left(  v\right)  .
\]
Furthermore since $\left\{  1\right\}  \oplus SU\left(  n-1\right)  $ is a
Poisson-Lie subgroup of $SU\left(  n\right)  $,
\[
\pi\left(  u^{\prime}\right)  =0\oplus\pi^{\left(  n-1\right)  }\left(
u^{\prime\prime}\right)  \in\left\{  0\right\}  \oplus\wedge^{2}%
T_{u^{\prime\prime}}SU\left(  n-1\right)  \subset\wedge^{2}T_{u^{\prime}%
}SU\left(  n\right)
\]
where $\pi^{\left(  n-1\right)  }$ is the standard multiplicative Poisson
structure on $SU\left(  n-1\right)  $. Note that
\[
\rho^{\left(  n-1\right)  }\left(  \phi_{3}\left(  v\right)  \right)  =\left(
D\psi\right)  _{u^{\prime\prime}}\left(  \pi^{\left(  n-1\right)  }\left(
u^{\prime\prime}\right)  \right)
\]
for the standard Poisson 2-tensor $\rho^{\left(  n-1\right)  }$ on
$\mathbb{S}^{2n-3}$. (Here we take $\mathbb{S}^{2n-3}=SU\left(  n-1\right)
/\left[  SU\left(  n-2\right)  \oplus\left\{  1\right\}  \right]  $.)

For
\[
u:=R_{\sigma_{c}}\left(  u^{\prime}\right)  =u^{\prime}\sigma_{c}\in SU\left(
n\right)  ,
\]
we have
\[
\phi_{1}\left(  u\right)  =u_{1}=\left(  u^{\prime}\sigma_{c}\right)
_{1}=v\in S_{c},
\]
and in $\wedge^{2}T_{v}S_{+}$,
\[
\left(  D\phi_{1}\right)  _{u}\left(  \pi_{\sigma_{c}}\left(  u\right)
\right)  =\left(  D\phi_{1}\right)  _{u}\left(  R_{\sigma_{c}}\left(
\pi\left(  u^{\prime}\right)  \right)  \right)  =\left(  D\phi_{1}\right)
_{u}\left(  \pi\left(  u^{\prime}\right)  \sigma_{c}\right)
\]%
\[
=\sqrt{1-c}\left(  D\psi\right)  _{u^{\prime}}\left(  \pi\left(  u^{\prime
}\right)  \right)  \in\wedge^{2}T_{v}S_{c}\subset\wedge^{2}T_{v}S_{+}%
\]
because the first columns of the component matrices in the 2-tensor
$\pi\left(  u^{\prime}\right)  =0\oplus\pi^{\left(  n-1\right)  }\left(
u^{\prime\prime}\right)  $ are all zero.

Note that
\[
\tau_{c}\left(  \left[  v\right]  \right)  =\tau_{c}\left(  \phi_{2}\left(
v\right)  \right)  =\tau_{c}\left(  \phi\left(  u\right)  \right)
\]%
\[
=\left(  D\phi\right)  _{u}\left(  \pi_{\sigma_{c}}\left(  u\right)  \right)
=\left(  D\phi_{2}\right)  _{\phi_{1}\left(  u\right)  }\left(  \left(
D\phi_{1}\right)  _{u}\left(  \pi_{\sigma_{c}}\left(  u\right)  \right)
\right)
\]
is a well-defines 2-tensor at $\left[  v\right]  \in\phi_{2}\left(
S_{c}\right)  \subset\mathbb{C}P^{n-1}$ and $\phi_{2}$ is a diffeomorphism on
$S_{+}$. So
\[
\pi^{\prime}:v\in S_{c}\mapsto\left(  D\phi_{1}\right)  _{u}\left(
\pi_{\sigma_{c}}\left(  u\right)  \right)  \in\wedge^{2}T_{v}S_{c}%
\]
is a well-defined Poisson 2-tensor on $S_{c}$ and $\phi_{2}\left(
S_{c}\right)  $ is a Poisson submanifold of $\left(  \mathbb{C}P^{n-1}%
,\tau_{c}\right)  $ that is Poisson isomorphic to $\left(  S_{c},\pi^{\prime
}\right)  $. Under the diffeomorphism $\phi_{3}:S_{c}\rightarrow
\mathbb{S}^{2n-3}$ identifying $v\in S_{c}$ with $\phi_{3}\left(  v\right)
\in\mathbb{S}^{2n-3}$, the 2-tensor $\pi^{\prime}\left(  v\right)  $ is
identified with
\[
\left(  D\phi_{3}\right)  _{v}\left(  \left(  D\phi_{1}\right)  _{u}\left(
\pi_{\sigma_{c}}\left(  u\right)  \right)  \right)  =\left(  D\phi_{3}\right)
_{v}\left(  \sqrt{1-c}\left(  D\psi\right)  _{u^{\prime}}\left(  \pi\left(
u^{\prime}\right)  \right)  \right)
\]%
\[
=\left(  D\phi_{3}\right)  _{v}\left(  \sqrt{1-c}\left(  0\oplus\left(
D\psi\right)  _{u^{\prime}}\left(  \pi^{\left(  n-1\right)  }\left(
u^{\prime\prime}\right)  \right)  \right)  \right)
\]%
\[
\left(  D\phi_{3}\right)  _{v}\left(  \sqrt{1-c}\left(  0\oplus\rho^{\left(
n-1\right)  }\left(  \phi_{3}\left(  v\right)  \right)  \right)  \right)
=\rho^{\left(  n-1\right)  }\left(  \phi_{3}\left(  v\right)  \right)
\in\wedge^{2}T_{\phi_{3}\left(  v\right)  }\mathbb{S}^{2n-3}.
\]
Thus $\left(  S_{c},\pi^{\prime}\right)  $ or $\left(  \phi_{2}\left(
S_{c}\right)  ,\tau_{c}\right)  $ is Poisson isomorphic to the standard
Poisson sphere $\left(  \mathbb{S}^{2n-3},\rho^{\left(  n-1\right)  }\right)
$.

When $n=2$, for $v\in S_{c}$ with $v_{2}\neq\sqrt{1-c}$, we cannot find a
$u^{\prime}\in SU\left(  2\right)  $ with the first column $u_{1}^{\prime
}=\left(  1,0\right)  $ and the last column $u_{2}^{\prime}=\sqrt{1-c}%
^{-1}\left(  0,v_{2}\right)  $. But for $v_{0}=\left(  \sqrt{c},\sqrt
{1-c}\right)  $, such a $u_{0}^{\prime}$ exists, namely, $u_{0}^{\prime}%
=I_{2}$ the $2\times2$ identity matrix, and the above argument essentially
works. More precisely, it is well known that $\pi\left(  u_{0}^{\prime
}\right)  =0$ since $u_{0}^{\prime}\in U\left(  1\right)  \subset SU\left(
2\right)  $, and hence for $u_{0}=u_{0}^{\prime}\sigma_{c}=\sigma_{c}$,
\[
\left(  D\phi_{1}\right)  _{u_{0}}\left(  \pi_{\sigma_{c}}\left(
u_{0}\right)  \right)  =\left(  D\phi_{1}\right)  _{u_{0}}\left(
R_{\sigma_{c}}\left(  \pi\left(  u_{0}^{\prime}\right)  \right)  \right)
=\left(  D\phi_{1}\right)  _{u_{0}}\left(  0\right)  =0.
\]
So
\[
\tau_{c}\left(  \left[  v_{0}\right]  \right)  =\left(  D\phi_{2}\right)
_{v_{0}}\left(  D\phi_{1}\right)  _{u_{0}}\left(  \pi_{\sigma_{c}}\left(
u_{0}\right)  \right)  =0.
\]
On the other hand, since $\tau_{c}$ on $\mathbb{C}P^{1}\approx\mathbb{S}^{2}$
is $SU\left(  2\right)  $-covariant and $U\left(  1\right)  \subset SU\left(
2\right)  $ consists of $0$-dimensional leaves, the action of any
\[
t=\left(
\begin{array}
[c]{cc}%
e^{i\theta} & 0\\
0 & e^{-i\theta}%
\end{array}
\right)  \in U\left(  1\right)
\]
on $\mathbb{C}P^{1}$ preserves the Poisson structure $\tau_{c}$. In
particular, $\tau_{c}\left(  \left[  tv_{0}\right]  \right)  =0$ for any $t\in
U\left(  1\right)  $. Since any $v=\left(  \sqrt{c},\sqrt{1-c}e^{i\theta
}\right)  \in S_{c}$ is equivalent to a $tv_{0}$ with $t\in U\left(  1\right)
$ under the diagonal $\mathbb{T}$-action, namely,
\[
\left[  v\right]  =\left[  \left(
\begin{array}
[c]{c}%
e^{-i\theta/2}\sqrt{c}\\
e^{-i\theta/2}\sqrt{1-c}e^{i\theta}%
\end{array}
\right)  \right]  =\left[  \left(
\begin{array}
[c]{c}%
e^{-i\theta/2}\sqrt{c}\\
e^{i\theta/2}\sqrt{1-c}%
\end{array}
\right)  \right]
\]%
\[
=\left[  \left(
\begin{array}
[c]{cc}%
e^{-i\theta/2} & 0\\
0 & e^{i\theta/2}%
\end{array}
\right)  \left(
\begin{array}
[c]{c}%
\sqrt{c}\\
\sqrt{1-c}%
\end{array}
\right)  \right]  =\left[  \left(
\begin{array}
[c]{cc}%
e^{-i\theta/2} & 0\\
0 & e^{i\theta/2}%
\end{array}
\right)  v_{0}\right]
\]
in $\mathbb{C}P^{1}$, we have $\tau_{c}\left(  \left[  v\right]  \right)  =0$
for all $v\in S_{c}$, i.e. $\tau_{c}=0$ on $\phi_{2}\left(  S_{c}\right)
\subset\mathbb{C}P^{1}$. Since $\phi_{2}\left(  S_{c}\right)  $ is
diffeomorphic to $S_{c}$ and hence to $\mathbb{S}^{1}$, we get the standard
(trivial) Poisson $\left(  \mathbb{S}^{1},\rho^{\left(  1\right)  }\right)  $
embedded in $\left(  \mathbb{C}P^{1},\tau_{c}\right)  $.
\end{proof}

It is easy to see that the difference of covariant Poisson 2-tensors is an
invariant 2-tensor, though not necessarily a Poisson 2-tensor. In the next
section, we show that any (left) $SU\left(  n\right)  $-invariant 2-tensor
$\tilde{\tau}$ on $\mathbb{C}P^{n-1}$, i.e.
\[
L_{u}\left(  \tilde{\tau}\left(  x\right)  \right)  =\tilde{\tau}\left(
ux\right)
\]
for all $u\in SU\left(  n\right)  $ and $x\in\mathbb{C}P^{n-1}$, is a constant
multiple of the canonical $SU\left(  n\right)  $-invariant symplectic 2-tensor
on $\mathbb{C}P^{n-1}$, and hence a Poisson 2-tensor on $\mathbb{C}P^{n-1}$.
So $\tau_{1}-\tau_{c}$ on $\mathbb{C}P^{n-1}$ for $c\in\left(  0,1\right)  $
is actually an $SU\left(  n\right)  $-invariant Poisson 2-tensor.

\section{Invariant 2-tensor on $\mathbb{S}^{2n-1}$}

In this section, we first classify the $SU\left(  n\right)  $-invariant
(contravariant) 2-tensor on $\mathbb{S}^{2n-1}$, and then we conclude that the
canonical $SU\left(  n\right)  $-invariant symplectic structure on
$\mathbb{C}P^{n-1}$ gives the only, up to a constant factor, $SU\left(
n\right)  $-invariant 2-tensor on $\mathbb{C}P^{n-1}$.

For each $p\in\mathbb{S}^{2n-1}$, we have $ip\in T_{p}\mathbb{S}^{2n-1}$, and
the orthogonal complement $E_{p}:=\left\{  p,ip\right\}  ^{\perp}\subset
T_{p}\mathbb{S}^{2n-1}$ is a complex subspace of $\mathbb{C}^{n}%
=T_{p}\mathbb{C}^{n}$ endowed with a canonical symplectic structure
$\tilde{\Omega}_{p}$ determined by the complex hermitian structure on
$\mathbb{C}^{n}$. Indeed $\left(  d\omega\right)  _{p}=\tilde{\Omega}_{p}$ on
$E_{p}$ for the unique 1-form $\omega$, the standard contact structure, on
$\mathbb{S}^{2n-1}$ such that $\omega_{p}\left(  ip\right)  =1$ and
$\omega_{p}\left(  E_{p}\right)  =\left\{  0\right\}  $ at each $p\in
\mathbb{S}^{2n-1}$. The contact manifold $(\mathbb{S}^{2n-1},\omega)$ with the
diagonal $\mathbb{T}$-action on $\mathbb{S}^{2n-1}$ is the standard
prequantization \cite{Kos,We:sm} of the canonical $SU\left(  n\right)
$-invariant symplectic structure on $\mathbb{C}P^{n-1}\cong\mathbb{S}%
^{2n-1}/\mathbb{T}$.

Since the vector fields $p\mapsto p$ and $p\mapsto ip$ on $\mathbb{S}^{2n-1}$
are invariant under the $U\left(  n\right)  $-action, so is the distribution
$p\mapsto E_{p}$ of tangent subspaces. Furthermore, since the $U\left(
n\right)  $-action preserves the complex hermitian structure on $\mathbb{C}%
^{n}$ (and on $E_{p}$), the field $p\mapsto\tilde{\Omega}_{p}$ of symplectic
forms on $\mathbb{S}^{2n-1}$ is also invariant under the $U\left(  n\right)
$-action. Thus the contravariant 2-tensor $\tilde{\pi}$ on $\mathbb{S}^{2n-1}$
uniquely determined by the form $\tilde{\Omega}$ on $E\subset TSU\left(
n\right)  $ is $U\left(  n\right)  $-invariant. Note that this 2-tensor
$\tilde{\pi}$ on $\mathbb{S}^{2n-1}$, invariant under the diagonal
$\mathbb{T}$-action, induces the cnanonical symplectic structure on
$\mathbb{C}P^{n-1}\cong\mathbb{S}^{2n-1}/\mathbb{T}$ determined by its complex
hermitian structure.

Given an $SU\left(  n\right)  $-invariant contravriant 2-tensors $\pi\neq0$ on
$\mathbb{S}^{2n-1}$ with $n\geq2$, we show that $\pi=\tilde{\pi}$ after a
suitable normalization if $n\neq3$ or if $\pi$ is $U\left(  n\right)
$-invariant. Through the standard Euclidean structure on $\mathbb{C}^{n}%
\cong\mathbb{R}^{2n}$, we identify the $SU\left(  n\right)  $-invariant
contravriant 2-tensors $\pi\neq0$ on $\mathbb{S}^{2n-1}$ with an $SU\left(
n\right)  $-invariant 2-forms $\Omega\neq0$ on $\mathbb{S}^{2n-1}$.

First we show that the tangent vector
\[
e_{1}^{\prime}:=ie_{1}\in T_{e_{1}}\mathbb{S}^{2n-1}=i\mathbb{R}%
\oplus\mathbb{C}^{n-1}%
\]
at $e_{1}\in\mathbb{S}^{2n-1}$ is in
\[
\ker\Omega_{p}:=\left\{  v\in T_{p}\mathbb{S}^{2n-1}:\Omega_{p}\left(
v,\cdot\right)  =0\right\}  .
\]
If not, then we can find an orthonormal set $\left\{  e_{i}^{\prime}\right\}
_{i=2}^{n-1}\cup\left\{  \eta_{i}^{\prime}\right\}  _{i=1}^{n}\subset
0\oplus\mathbb{C}^{n-1}$ such that $\Omega_{e_{1}}\left(  e_{i}^{\prime}%
,\eta_{j}^{\prime}\right)  =\delta_{ij}a_{ii}$ and $\Omega_{e_{1}}\left(
e_{i}^{\prime},e_{j}^{\prime}\right)  =\Omega_{e_{1}}\left(  \eta_{i}^{\prime
},\eta_{j}^{\prime}\right)  =0$ with $a_{ii}\in\mathbb{R}$ and $a_{11}\neq0$.
Now since $\Omega$ is $SU\left(  n\right)  $-invariant, we have
\[
\Omega_{e_{1}}\left(  e_{1}^{\prime},u\left(  \eta_{1}^{\prime}\right)
\right)  =\Omega_{u\left(  e_{1}\right)  }\left(  u\left(  e_{1}^{\prime
}\right)  ,u\left(  \eta_{1}^{\prime}\right)  \right)  =\Omega_{e_{1}}\left(
e_{1}^{\prime},\eta_{1}^{\prime}\right)  =a_{11}%
\]
for any $u\in\left\{  1\right\}  \oplus SU\left(  n-1\right)  \subset
SU\left(  n\right)  $. This cannot be true, since by a suitable choice of $u$,
$u\left(  \eta_{1}^{\prime}\right)  $ can be any unit vector in $0\oplus
\mathbb{C}^{n-1}$, for example, $\eta_{n}^{\prime}$. Thus $e_{1}^{\prime
}=ie_{1}\in\ker\Omega_{p}$.

Now with respect to the standard orthonormal $\mathbb{R}$-linear basis of
\[
i\mathbb{R}\oplus\mathbb{R}^{n-1}\oplus\mathbb{R}^{n-1}\cong i\mathbb{R}%
\oplus\mathbb{C}^{n-1}=T_{e_{1}}\mathbb{S}^{2n-1},
\]
the 2-form $\Omega_{e_{1}}$ can be represented by a block diagonal matrix
$0\oplus B$ where $B\in M_{2\left(  n-1\right)  }\left(  \mathbb{R}\right)  $
is a skew symmetric matrix. The $SU\left(  n\right)  $-invariance of
$\Omega\neq0$ implies that $\Omega_{e_{1}}\neq0$ and
\[
uBu^{-1}=uBu^{t}=B\neq0,
\]
or $uB=Bu$, for any $u\in SU\left(  n-1\right)  \subset O_{2n-2}\left(
\mathbb{R}\right)  $ since $1\oplus u\in SU\left(  n\right)  $ and $\left(
1\oplus u\right)  \left(  e_{1}\right)  =e_{1}$. (If $\Omega$ is $U\left(
n\right)  $-invariant, then $uB=Bu$ for any $u\in U\left(  n-1\right)  $ since
$1\oplus u\in U\left(  n\right)  $.)

We claim that $B$ must be conformal, i.e. $\left|  \left|  B\left(  v\right)
\right|  \right|  =\left|  \left|  B\right|  \right|  \left|  \left|
v\right|  \right|  $ for all $v\in\mathbb{R}^{2n-2}$ where $\left|  \left|
B\right|  \right|  :=\sup_{\left|  \left|  v\right|  \right|  =1}\left|
\left|  B\left(  v\right)  \right|  \right|  >0$. Let $w$ be a unit vector
with $\left|  \left|  B\left(  w\right)  \right|  \right|  =\left|  \left|
B\right|  \right|  $. Since $SU\left(  n-1\right)  $ acts on $\mathbb{S}%
^{2n-3}\subset\mathbb{R}^{2n-2}$ transitively, for any unit vector
$v\in\mathbb{R}^{2n-2}$, we can find $u\in SU\left(  n-1\right)  $ with
$u^{-1}\left(  v\right)  =w$, and hence
\[
\left|  \left|  B\left(  v\right)  \right|  \right|  =\left|  \left|
uBu^{-1}\left(  v\right)  \right|  \right|  =\left|  \left|  u\left(  B\left(
w\right)  \right)  \right|  \right|  =\left|  \left|  B\left(  w\right)
\right|  \right|  =\left|  \left|  B\right|  \right|  .
\]
Thus $B/\left|  \left|  B\right|  \right|  $ is a skew-symmetric isometry on
$\mathbb{R}^{2n-2}$ and so $B/\left|  \left|  B\right|  \right|  \in
O_{2n-2}\left(  \mathbb{R}\right)  $.

If $n=2$, then any skew symmetric $0\neq B/\left|  \left|  B\right|  \right|
\in O_{2}\left(  \mathbb{R}\right)  $ determines the same 2-form
$\Omega_{e_{1}}$ on $0\oplus\mathbb{R}^{2}$ and hence on $i\mathbb{R}%
\oplus\mathbb{R}^{2}$, up to a constant multiple. So $\Omega=\tilde{\Omega}$
after normalized.

If $n\geq4$, then the commutativity of $\mathbb{T}^{n-2}\subset SU\left(
n-1\right)  $ with $B$ implies that $B$ is complex linear on $\mathbb{R}%
^{2n-2}=\mathbb{C}^{n-1}$ and so $B/\left|  \left|  B\right|  \right|  \in
U\left(  n-1\right)  $. In fact, since for any $1\leq j\neq k\leq n-1$,
$t_{jk\theta}B=Bt_{jk\theta}$ for all $\theta\in\mathbb{R}$ implies that
$B_{jj},B_{kk}\in\mathbb{C}$ and $B_{kl}=0$ for any $j\neq l\neq k$, where
$B=\left(  B_{jk}\right)  _{1\leq j,k\leq n-1}$ with $B_{jk}\in M_{2}\left(
\mathbb{R}\right)  $, and
\[
t_{jk\theta}:=e^{i\theta}e_{jj}+e^{-i\theta}e_{kk}+\sum_{\substack{1\leq l\leq
n-1 \\j\neq l\neq k}}e_{ll}\in\mathbb{T}^{n-2}\subset SU\left(  n-1\right)  .
\]
It is well known that only scalar matrices in $M_{n-1}\left(  \mathbb{C}%
\right)  $ commute with $SU\left(  n-1\right)  $, so we get $B/\left|  \left|
B\right|  \right|  \in\mathbb{T}$ with $-B/\left|  \left|  B\right|  \right|
=\left(  B/\left|  \left|  B\right|  \right|  \right)  ^{\ast}=\left(
B/\left|  \left|  B\right|  \right|  \right)  ^{-1}$, i.e. $\left(  B/\left|
\left|  B\right|  \right|  \right)  ^{2}=-1$ or $B=\pm i\left|  \left|
B\right|  \right|  $. Thus
\[
\Omega_{e_{1}}=\pm\left|  \left|  B\right|  \right|  \tilde{\Omega}_{e_{1}}%
\]
a (real) constant multiple of the standard symplectic form. Hence we get
$\pi=\tilde{\pi}$ after a suitable normalization.

If $\Omega$ is $U\left(  n\right)  $-invariant, then the commutativity of
$\mathbb{T}^{n-1}\subset U\left(  n-1\right)  $ with $B$ implies that $B$ is
complex linear and hence $B/\left|  \left|  B\right|  \right|  \in U\left(
n-1\right)  $ and as above, $\Omega=\pm\left|  \left|  B\right|  \right|
\tilde{\Omega}$. In fact, $t_{k\theta}^{\prime}B=Bt_{k\theta}^{\prime}$ for
all $\theta\in\mathbb{R}$ implies that $B_{kk}\in\mathbb{C}$ and $B_{kl}=0$
for any $l\neq k$, where
\[
t_{k\theta}^{\prime}:=e^{i\theta}e_{kk}+\sum_{\substack{1\leq l\leq n-1
\\l\neq k}}e_{ll}\in\mathbb{T}^{n-1}\subset U\left(  n-1\right)  .
\]

We observe that the quotient map $\phi:\mathbb{S}^{2n-1}\rightarrow
\mathbb{C}P^{n-1}$ and its differential $D\phi:T\mathbb{S}^{2n-1}\rightarrow
T\mathbb{C}P^{n-1}$ are $U\left(  n\right)  $-equivariant since the diagonal
$\mathbb{T}$-action commutes with the $U\left(  n\right)  $-action.
Furthermore, the restriction
\[
\left(  D\phi\right)  |_{E}:E\rightarrow T\mathbb{C}P^{n-1}%
\]
of $D\phi$ to the $U\left(  n\right)  $-equivariant subbundle $E$ defined
above is a bundle isomorphism. So any $SU\left(  n\right)  $-invariant (and
hence $U\left(  n\right)  $-invariant) 2-tensor $\tau\in\Gamma\left(
\wedge^{2}T\mathbb{C}P^{n-1}\right)  $ on $\mathbb{C}P^{n-1}$ can be `pulled
back' to an $U\left(  n\right)  $-invariant 2-tensor
\[
\pi=\left(  D\phi\right)  |_{E}^{-1}\left(  \tau\right)  \in\Gamma\left(
\wedge^{2}E\right)  \subset\Gamma\left(  \wedge^{2}T\mathbb{S}^{2n-1}\right)
\]
on $\mathbb{S}^{2n-1}$ which must be, up to a constant factor, equal to
$\tilde{\pi}$. Thus $\tau=\tilde{\tau}:=\left(  D\phi\right)  \left(
\tilde{\pi}\right)  $ which is the standard symplectic 2-tensor on
$\mathbb{C}P^{n-1} $.


\begin{thebibliography}{99}
\bibitem[ DaSo]{DaSo}P. Dazord and D. Sondaz, \emph{Groupes de Poisson
affines}, in `Symplectic Geometry, Groupoids, and Integrable Systems', P.
Dazord and A. Weinstein (Eds.), Springer-Verlag, 1991.

\bibitem[ DiNo] {DiNo}M. S. Dijkhuizen and M. Noumi, \emph{A family of quantum
projective spaces and related }$q$\emph{-hypergeometric orthogonal
polynomials}, preprint.

\bibitem[ Dr] {D}V. G. Drinfeld, \emph{Quantum groups}, Proc. I.C.M. Berkeley
1986, Vol. 1, 789-820, Amer. Math. Soc., Providence, 1987.

\bibitem[ KhRaRu] {KRR}S. Khoroshkin, A. Radul, and V. Rubtsov, \emph{A family
of Poisson structures on hermitian symmetric spaces}, Comm. Math. Phys. 152
(1993), 299-315.

\bibitem[ KoVa] {KoVa}L. I. Korogodsky and L. L. Vaksman, \emph{Quantum }%
$G$\emph{-spaces and Heisenberg algebra}, in ``Quantum Groups'', P. P. Kulish
(Ed.), Lecture Notes in Mathematics 1510, Springer-Verlag, Berlin, 1992, pp. 56-66.

\bibitem[ Kos] {Kos}B. Kostant, \emph{Quantization and unitary
representations. I: Prequantization}, in ``Lectures in Modern Analysis and
Applications. III'', E. T. Taam (Ed.), Lecture Notes in Mathematics 170,
Springer-Verlag, Berlin and New York, 1970, pp. 87-208.

\bibitem[ Lu] {Lu:maps}J. H. Lu, \emph{Multiplicative and affine Poisson
structures on Lie groups}, Ph. D. thesis, Univ. of California, Berkeley, 1990.

\bibitem[ LuWe1] {LuWe:plg}J. H. Lu and A. Weinstein, \emph{Poisson Lie
groups, dressing transformations and Bruhat decompositions}, J. Diff. Geom. 31
(1990), 501-526.

\bibitem[ LuWe2] {LuWe:ps}\_\_\_\_\_\_, \emph{Classification of }%
$SU(2)$\emph{-covariant Poisson structures on }$\mathbb{S}^{2}$, Comm. Math.
Phys. 135 (1991), 229-231.

\bibitem[ Po] {Po}P. Podles, \emph{Quantum spheres}, Letters Math. Phys. 14
(1987), 193-202.

\bibitem[ Re] {Re}J. Renault, A Groupoid Approach to C*-algebras, Lecture
Notes in Mathematics, Vol. 793, Springer-Verlag, New York, 1980.

\bibitem[ RTF] {FRT}N. Yu. Reshetikhin, L. A. Takhtadzhyan, and L. D. Faddeev,
\emph{Quantization of Lie groups and Lie algebras}, Leningrad Math. J. 1
(1990), 193-225.

\bibitem[ Ri] {Ri:ncq}M. A. Rieffel, \emph{Non-compact quantum groups
associated with abelian subgroups}, Comm. Math. Phys. 171 (1995), 181-201.

\bibitem[ Sh1] {Sh:qp}A. J. L. Sheu, \emph{Quantization of the Poisson SU(2)
and its Poisson homogeneous space -- the 2-sphere}, Comm. Math. Phys. 135
(1991), 217-232.

\bibitem[ Sh2] {Sh:cqg}\_\_\_\_\_\_, \emph{Compact quantum groups and groupoid
C*-algebras}, J. Func. Anal. 144 (1997), 371-393.

\bibitem[ Sh3] {Sh:qcp}\_\_\_\_\_, \emph{Groupoid approach to quantum
projective spaces}, preprint.

\bibitem[ So] {So}Ya. S. Soibelman, \emph{The algebra of functions on a
compact quantum group, and its representations}, Algebra Analiz. 2 (1990),
190-221. (Leningrad Math. J., 2 (1991), 161-178.)

\bibitem[ VaSo] {VaSo:af}L. L. Vaksman and Ya. S. Soibelman, \emph{The algebra
of functions on the quantum group }$SU(n+1)$\emph{, and odd-dimensional
quantum spheres}, Leningrad Math. J. 2 (1991), 1023-1042.

\bibitem[ We1] {We:sm}A. Weinstein, \emph{Lectures on Symplectic Manifolds},
CBMS Regional Conference Series in Mathematics, No. 29, Amer. Math. Soc.,
Providence, 1977.

\bibitem[ We2] {We}A. Weinstein, \emph{The local structure of Poisson
manifolds}, J. Diff. Geom. 18 (1983), 523-557.

\bibitem[ We3] {We:aps}A. Weinstein, \emph{Affine Poisson structures}, preprint.

\bibitem[ Wo1] {Wo:ts}S. L. Woronowicz, \emph{Twisted $SU(2)$ group: an
example of a non-commutative differential calculus}, Publ. RIMS. 23 (1987), 117-181.

\bibitem[ Wo2] {Wo:cm}\_\_\_\_\_\_, \emph{Compact matrix pseudogroups}, Comm.
Math. Phys. 111 (1987), 613-665.
\end{thebibliography}
\end{document}